\documentclass[12pt]{amsart}
\usepackage{amsfonts,amssymb}
\title[Bounds on Multiplicities of Laplace Operator Eigenvalues]{Bounds on Multiplicities 
of Laplace-Beltrami Operator Eigenvalues on the Real Projective Plane}
\author[A.~S.~Berdnikov]{Aleksandr~S.~Berdnikov}
\address[A.~S.~Berdnikov]{Department of Mathematics, Massachusetts Institute of Technology, 
77 Massachusetts Ave, Cambridge, MA 02139, USA
\newline {\em and} 
\newline Independent University of Moscow, Bolshoy Vlasyevskiy per.~11, 119002, Moscow, Russia
\newline {\em and}
\newline National Research University Higher School of Economics,
Vavilova Str.~7, 117312, Moscow, Russia}
\email{\tt aberdnik@mit.edu}
\author[N.~S.~Nadirashvili]{Nikolai~S.~Nadirashvili}
\address[N.~S.~Nadirashvili]{Aix Marseille Universit\'e, CNRS,
I2M UMR 7353 --- Centre de Math\'e\-ma\-ti\-ques et Informatique, 13453, Marseille, France}
\email{nikolay.nadirashvili@univ-amu.fr}
\author[A.~V.~Penskoi]{Alexei V. Penskoi}\thanks{The work of the third author
was supported by Russian Science Foundation grant no.~16-11-10260
at Moscow State University}
\address[Alexei V. Penskoi]{Department of Higher Geometry and Topology, 
Faculty of Mathematics and Mechanics, Moscow State University,
Leninskie Gory, GSP-1, 119991, Moscow, Russia \newline {\em and}
\newline Faculty of Mathematics,
National Research University Higher School of Economics,
6 Usacheva Str., 119048, Moscow, Russia \newline {\em and}
\newline Laboratoire J.-V.Poncelet (UMI 2615), Bolshoy Vlasyevskiy 
Pereulok 11, 119002, Moscow, Russia}
\email[corresponding author]{penskoi@mccme.ru}
\date{}
\DeclareMathOperator{\chr}{chr}
\newtheorem{definition}{Definition}
\newtheorem{lema}{Lemma}
\newtheorem{prop}{Proposition}
\newtheorem{teo}{Theorem}
\begin{document}
\begin{abstract}
The known upper bounds for the multiplicities
of the Laplace-Beltrami operator eigenvalues
on the real projective plane
are improved for the eigenvalues with even indexes.
Upper bounds for Dirichlet, Neumann and Steklov
eigenvalues on the real projective plane with
holes are also provided.
\end{abstract}
\maketitle

\section*{Introduction}

Consider a closed surface $\Sigma$ with 
a Riemannian metric $g$. 
Let
$$
0=\lambda_0<\lambda_1\leqslant\lambda_2\leqslant\dots
$$
be Laplace-Beltrami operator eigenvalues (counting with multiplicities).
Let $m(\Sigma,g,\lambda_i)$ denote the multiplicity of eigenvalue $\lambda_i$, 
i.e. the dimension of the eigenspace corresponding to $\lambda_i$.

Finding upper bounds for $m(\Sigma,g,\lambda_i)$ is one of classical
problems of Spectral Geometry. These bounds are interesting
from several points of view. For example, they play an important role
in the study of metrics extremal for the Laplace-Beltrami
eigenvalues, see e.g.~\cite{Penskoi2013}.

In the present paper we consider the case of $\Sigma$ being the real 
projective plane $\mathbb{RP}^2.$
The best known bound for the real projective plane
\begin{equation}
\label{c}
m(\mathbb{R}P^2,g,\lambda_i)\leqslant 2i+3
\end{equation}
was proven in the paper~\cite{closed}, see Theorem~\ref{Nadirashvili87}
below. Recently this bound was improved in the paper~\cite{recent}
for $i=2,$
$$
m(\mathbb{R}P^2,g,\lambda_2)\leqslant 6,
$$
see Theorem~\ref{Nadirashvili-Penskoi} below.

The first goal of the present paper is to
improve the bounds~\eqref{c} for all
even $i$ and prove the following theorem.

\begin{teo}\label{m1}
Let $l$ be a positive integer, then
$$m(\mathbb{R}P^2,g,\lambda_{2l})\leqslant 4l+1.$$
\end{teo}

We also consider the case of a surface $\Sigma$ with a boundary 
$\partial \Sigma$. In this case, let 
$$
0\leqslant \lambda_0 \leqslant \lambda_1 \leqslant \dots
$$
denotes the spectrum of the Laplace-Beltrami operator, given 
either Dirichlet or Neumann boundary condition on each connected 
component of the boundary $\partial \Sigma$. 

We also consider Steklov eigenvalue problem on $(\Sigma,\partial \Sigma).$
Let $\rho$ be a bounded non-negative function on $\partial \Sigma$ and 
$\sigma$ a real number. Then a function $u$ on $\Sigma$ is called 
a Steklov eigenfunction with eigenvalue $\sigma$ if 
$$
\left\{\begin{array}{l}
\Delta u=0 \mbox{ in $\Sigma$,}\\
\dfrac{\partial u}{\partial n}=\sigma\rho u \mbox{ on $\partial \Sigma$,}
\end{array}\right.
$$
where $n$ denotes the unit outer normal on $\partial \Sigma$. 

Let us consider the case where $\rho$ is the density of an absolutely 
continuous Radon measure $s$ on $\partial \Sigma$. Then the spectrum of 
Steklov problem is non-negative and discrete~\cite{steklov}, 
so we denote it by
$$
0= \sigma_0 \leqslant \sigma_1 \leqslant \dots
$$
Let us denote the multiplicity corresponding to the Steklov eigenvalue 
$\sigma_i$ by $m(\Sigma,g,s,\sigma_i)$.

Let now $\mathbb{RP}^2_h$ denote $\mathbb{RP}^2\setminus (\cup D^2_i)$, 
i.e. the real projective plane with a positive number of holes. 
Our second result then is the following theorem. 

\begin{teo}
\label{m2}
Let $l$ be a positive integer and $s$ be an absolutely continuous Radon measure on $\partial \mathbb{R}P^2_h$ with bounded density, then
\begin{equation}\label{bound-DN}
m(\mathbb{R}P^2_h,g,\lambda_{2l})\leqslant 4l+2,
\end{equation}
\begin{equation}\label{bound-steklov}
m(\mathbb{R}P^2_h,g,s,\sigma_{2l})\leqslant 4l+2.
\end{equation}
\end{teo}

The proofs follow the technique from the paper~\cite{closed} but uses a 
more refined topological argument at the last step.

Let us recall some already known results about upper bounds on 
multiplicities of eigenvalues.

\begin{teo}[Cheng~\cite{cheng}]
Let $\Sigma$ be an oriented surface of genus $\gamma$. Then
for any metric $g$ one has 
$$
m(\Sigma,g,\lambda_i) \leqslant \frac{(2\gamma+i+1)(2\gamma+i+2)}{2}.
$$
\end{teo}

\begin{teo}[Besson~\cite{besso}]
Let $\Sigma$ be an oriented surface of genus $\gamma$. Then
for any metric $g$ one has
$$
m(\Sigma,g,\lambda_i) \leqslant 4\gamma+2i+1.
$$
Let $\Sigma$ be a non-orientable surface of Euler characteristic 
$\chi(\Sigma)$. Then
for any metric $g$ one has 
$$
m(\Sigma,g,\lambda_i)\leqslant 4(1-\chi(\Sigma))+4i+2.
$$
\end{teo}

\begin{teo}[Nadirashvili~\cite{closed}]\label{Nadirashvili87}
For any metric $g$  on the sphere $\mathbb{S}^2$, the real projective plane 
$\mathbb{RP}^2$, the torus $\mathbb{T}^2$, or the Klein bottle $\mathbb{KL}$
the following inequalities hold,
$$
m(\mathbb{S}^2,g,\lambda_i) \leqslant 2i+1,$$
$$m(\mathbb{RP}^2,g,\lambda_i) \leqslant 2i+3,$$
$$m(\mathbb{T}^2,g,\lambda_i) \leqslant 2i+4,$$
$$m(\mathbb{KL},g,\lambda_i) \leqslant 2i+3.
$$
For any other surface $\Sigma,$ i.e. for a surface $\Sigma$ with 
$\chi(\Sigma)<0,$ with any metric $g$ the following inequality holds,
$$
m(\Sigma,g,\lambda_i) \leqslant 2i-2\chi(\Sigma)+3.
$$
\end{teo}

\begin{teo}[M. Hoffmann-Ostenhof, T. Hoffmann-Ostenhof, and 
Nadirashvili~\cite{hhn}]
Let $\Sigma$ be a closed surface of genus 0. Then for any metric $g$ the inequality
$$
m(\Sigma,g,\lambda_i) \leqslant 2i - 3
$$
holds for $i\geqslant 3$.
\end{teo}

\begin{teo}[Nadirashvili and Penskoi~\cite{recent}]\label{Nadirashvili-Penskoi}
The following upper bound for the multiplicity of the second 
eigenvalue of the Laplace-Beltrami operator on the projective plane holds
for any metric $g,$
$$
m(\mathbb{R}P^2,g,\lambda_2)\leqslant 6.
$$
\end{teo}

For a surface $\Sigma$ with boundary $\partial \Sigma$ let us
denote by $\overline{\Sigma}$ the closed (topological) manifold, 
obtained by collapsing each connected component of $\partial \Sigma$ 
into a point.

\begin{teo}[Karpukhin, Kokarev, and Polterovich~\cite{polt}]
Let $(M,g)$ be a compact Riemannian surface with a non-empty boundary, 
and $\mu$ be an absolutely continuous Radon measure on $\partial M$ 
whose density is bounded. Then the multiplicity $m(M,g,\mu,\sigma_k)$ 
of a Steklov eigenvalue $\sigma_k(g,\mu)$ satisfies the inequalities
\begin{equation}\label{kkp1}
m(M,g,\mu,\sigma_k) \leqslant 2(2- \chi (\overline{M})) +2k+1,
\end{equation}
\begin{equation}\label{kkp2}
m(M,g,\mu,\sigma_k) \leqslant 2(2- \chi (M)) +k,
\end{equation}
for all $k = 1,2\dots$. 

Besides, the latter inequality is strict for an even $k$.
\end{teo}

\begin{teo}[Jammes~\cite{jammes}]
Let $\Sigma$ be a compact surface with boundary, then for 
$k\geqslant 1$ the following inequality holds,
\begin{equation}\label{bound-jammes}
m(\Sigma,g,\mu,\sigma_k)\leqslant k-2\chi(\Sigma)+3.
\end{equation}
\end{teo}

\begin{teo}[T. Hoffmann-Ostenhof, Michor, and Nadirashvili~\cite{orig}]
Let $k\geqslant 3$. Then the multiplicity of the $k$-th eigenvalue 
$\lambda_k$ for the Dirichlet problem on a planar simply-connected 
domain $D$ satisfies the inequality
$$m(D,\lambda_k) \leqslant 2k-1.$$
\end{teo}

\begin{teo}[Berdnikov~\cite{berd}]
Let $M$ be the surface with $\chi(\overline{M})<0$. Then 
$$m(M,g,\lambda_k)\leqslant 2k-2\chi (\overline{M})+3,$$
$$m(M,g,vol_{\partial M},\sigma_k)\leqslant 2k-2\chi (\overline{M})+3.
$$
\end{teo}

Let us also recall results concerning the relation between
bounds on multiplicities of eigenvalues and the chromatic number.

\begin{definition} A chromatic number $\chr(\Sigma)$ of a surface $\Sigma$
is the maximal $n$ such that one can embed in $\Sigma$ the complete graph 
$K_n$ on $n$ vertices.
\end{definition}

Let us consider a Schr\"odinger operator $H=\Delta+V.$
Let $\bar{m}_1(\Sigma)$ denote the supremum over all $H$
of the multiplicitiy of the eigenvalue $\lambda_1$ on $\Sigma.$

\begin{teo}[Colin de Verdi\`ere~\cite{CdV1987}]
For any surface $\Sigma,$ one has $\bar{m}_1(\Sigma)\geqslant \chr(\Sigma)-1.$

For the four surfaces $\Sigma$ with $\chi(S)\geqslant0,$ namely 
$\mathbb{S}^2,$ $\mathbb{RP}^2,$ $\mathbb{T}^2,$
$\mathbb{KL},$ one has $\bar{m}_1(\Sigma)=\chr(\Sigma)-1.$
\end{teo}

This theorem leads to a natural conjecture that for any surface $\Sigma$
one has $\bar{m}_1(\Sigma)=\chr(\Sigma)-1.$ This conjecture was proved
by S\'evennec for some surfaces $\Sigma$ using the following upper bound
for $\bar{m}_1(\Sigma).$

\begin{teo}[S\'evennec~\cite{Sevennec2002}] If $\chi(\Sigma)<0$ then
$\bar{m}_1(\Sigma)\leqslant5-\chi(\Sigma).$ 

It follows that the above conjecture holds for
all surfaces $\Sigma$ with $\chi(\Sigma)\geqslant-3,$
the four new cases being $\mathbb{T}^2\#\mathbb{T}^2$ and 
$\#n\mathbb{RP}^2,$ $n =3,4,5.$
\end{teo}

Since $m(\Sigma,g,\lambda_1)\leqslant\bar{m}_1(\Sigma),$ these results provide
interesting upper bounds on $m(\Sigma,g,\lambda_1).$

Let us now compare Theorem~\ref{m2} with the above mentioned results.
The bound given by formula~\eqref{bound-DN}
seems to be a first bound of this kind for the projective plane with
holes.

The comparison of the bound given in formula~\eqref{bound-steklov}
with bounds~\eqref{kkp1}, \eqref{kkp2} and \eqref{bound-jammes}
gives a complicated answer. Bounds~\eqref{kkp1}, \eqref{kkp2}
are linear in $k,$ bound~\eqref{kkp1} has a better constant,
bound~\eqref{kkp2} has a better asymptotic. Bound~\eqref{bound-steklov}
improves bound~\eqref{kkp1} for $\mathbb{RP}^2$ by $1.$
Bound~\eqref{bound-jammes} improves bound~\eqref{kkp2}
by $1.$

The third author is very grateful for the
the Institut de Ma\-th\'e\-ma\-ti\-ques de Marseille (I2M, UMR~7373)
for the hospitality.
The third author is a Young Russian Mathematics award winner 
and would like to thank its sponsors and jury.
The third author is a Simons research-professorship award winner and 
would like to express his deep gratitude to the jury and to 
the Simons foundation.
The third author is very indebted to Pavel Winternitz
for fruitful discussions.

\section{Preliminaries}

The standard relation between the original problem and the topology of the 
surface is well-known. Recall that for an eigenfunction $f$ of the Laplacian, 
the connected components of $f^{-1}(\mathbb{R}\setminus \{ 0\})$ are 
called {\itshape nodal domains} of $f$ and the set $f^{-1}(0)$ is 
called {\itshape a nodal graph} $\mathcal{N}(f)$ of $f$. We will denote the 
number of nodal domains of a (nodal) graph $\mathcal{N}$ by $\mu(\mathcal{N})$. 
We will also use notation $\overline{\mathcal{N}}(f)$ 
for $\overline{\mathcal{N}(f)\cap Int(\Sigma)}\subset \overline{\Sigma}$. 
The term ``graph'' is justified by the following theorem.

\begin{teo}[Bers theorem, \cite{bers}] For a Laplacain 
eigenfunction $f$ and $x_0\in M$ there exists an integer $n\geqslant 0$ 
and polar coordinates $(r,\theta)$ centered at $x_0$ such 
that the following formula holds
$$
f(x)=r^n(\sin (n\theta))+O(r^{n+1}).
$$
\end{teo}

The Bers theorem implies that the nodal graph is, in fact, an embedded graph. 
It also implies that each eigenfunction has a well-defined order $ord(f,x):=n$ 
of vanishing at any point $x_0\in M$. Note that the quotient of eigenfunctions
$f$ such that $ord(f,x)=n$ by those $g$ with $ord(g,x)=n+1$ is contained in 
the span of $r^n(\sin (n\theta))$ and $r^n(\cos (n\theta))$ and hence at 
most 2-dimensional (if $n=0$ then $\sin (n\theta)=\cos (n\theta)$ and the 
quotient is at most 1-dimensional). That implies that high-dimensional 
eigenspaces contain functions with high vanishing order at a given point.

\begin{prop}
\label{tec1}
Let $U$ be an eigenspace of Laplacian with $\dim (U)\geqslant n$ and 
$x$ be a point in $M$. Then the dimension of the linear subspace of 
functions $f\in U$ such that $\mbox{\upshape ord}(f,x)\geqslant k>0$ 
(or, equivalently, $x$ is a vertex of $\mathcal{N}(f)$ of 
$\deg \geqslant 2k$) is at least $1+n-2k$.
\end{prop}

The following theorem is crucial for the considered approach, it contains 
an upper bound for the number of nodal domains.

\begin{teo}[Courant nodal domain theorem, \cite{c}, \cite{courant}, 
\cite{courantsteklov}]
Let $\Sigma$ be a smooth manifold with smooth boundary $\partial \Sigma$ 
(possibly empty) and $s$ be an absolutely continuous Radon measure on 
$\partial M$ with bounded density. Consider either the Laplacian 
eigenvalue problem on $M$ with the Dirichlet or Neumann boundary condition 
on each connected component of $\partial \Sigma$, or, if 
$\partial \Sigma \neq \varnothing$, the Steklov eigenvalue problem for the 
measure $s$. Then for each non-zero function $f$ in the $i$-th eigenspace 
$U_i$ of the considered spectral problem the number $\mu(\mathcal{N}(f))$ 
of its nodal domains is not greater than $i+1$.
\end{teo}

The Theorems~\ref{m1} and~\ref{m2} follow immediately  
from the Courant nodal domain theorem and the following theorem, which we 
prove later in this paper.

\begin{teo}
\label{B}
Let $\lambda$ be a real number and $U$ be a linear space of functions on 
the closed surface $\Sigma=\mathbb{R}P^2$ or a surface with the 
boundary $\Sigma=\mathbb{R}P^2_h$. Suppose that every $f\in U$ is a 
Laplacian eigenfunction $\Delta f=\lambda f$. Suppose 
$\sup\limits_{f\in U} (\mu (\mathcal{N} (f)))$ is not greater than an 
odd number~$n=2l+1.$ Let $i=n-1$. Then $\dim (U)\leqslant 2n=2i+2$. 

Moreover, if $\Sigma=\mathbb{R}P^2$, i. e. $\partial \Sigma=\varnothing$, 
then this inequality is strict, i.e. $\dim (U)<2n=2i+2$.
\end{teo}

We will use a notion of a star fibration, which we explain now. Bers theorem 
tells us that the nodal graph of the eigenfunction $f$ is diffeomorphic near 
$x_0$ to $2n$ rays in $\mathbb{R}^2$ emitting from $0$ at equal angles between 
the adjacent lines. That makes it natural to consider a star fibration 
$E_M(2n)$ over $Int(M)$ introduced in paper~\cite{berd}. 
It can be defined as 
follows. Recall that the spherisation $S(E)$ of a vector bundle $E$ is the 
fiber bundle $(E\setminus \underline{0})/\mathbb{R}_{>0}$. In the case when 
$E$ carries a positive-definite metric (i.e. tangent bundle of a
 Riemannian manifold), the spherisation $S(E)$ is isomorphic to the fiber 
bundle of the unit spheres in $E$. Now, consider the subset 
$E'_M(2n)\subset \prod^{2n} S(TM)$ such that the point 
$((x_1,v_1),\dots, (x_{2n},v_{2n}))$ belongs to $E'_M(2n)$ iff all $x_i$ 
are equal and $v_i$ are representing equidistant rays in $T_xM$. Define 
$E_M(2n)$ as a quotient of $E'_M(2n)$ by the natural action of the 
permutation group $S_n$. The fiber $F_x(2n)$ at a point $x$ of the 
fibration $E_M(2n)$ consists then of all $2n$-stars in the tangent space 
$T_xM$, i.e. configurations of $2n$ rays (or equally $n$ lines) in $T_xM$ 
with equal angles between adjacent lines. 

In these terms the Bers theorem states that if $ord(f,x)=n$ then the nodal 
graph $\mathcal{N}(f)$ defines an element of $F_x(2n)$ which we denote 
by $s(\mathcal{N}(f),x)$.

\section{First bound}

We start by proving Lemma 1 from paper~\cite{berd} for our particular case $\Sigma=\mathbb{RP}^2$ or $\Sigma=\mathbb{RP}^2_h$.

\begin{lema}
\label{cases}
Let $\lambda \in \mathbb{R}$ be a real number, $U$ be a linear space of 
functions on $\Sigma\in \{\mathbb{RP}^2,\mathbb{RP}^2_h\}$ satisfying 
$\Delta (f)=\lambda f$ and let 
$n=\sup\limits_{f\in U}  \mu (\mathcal{N} (f))$. 
Suppose that \mbox{$\dim (U)\geqslant 2n$}. For each $x\in Int(\Sigma)$ 
consider the set $U_n(x)\subset U$ consisting of eigenfunctions $f_x$ 
whose nodal graph $\mathcal{N}_f$ contains a vertex $x$ of degree at 
least $2n$.  

Then $\dim (U_n(x))\geqslant 1$. Moreover,
\begin{itemize}
\item any $f_{x}$ has a nodal graph with a unique vertex in $\overline{\Sigma}$;
\item $ord(f_x,x)=n$;
\item faces of $\overline{\mathcal{N}}(f_{x})$ are homeomorphic to $D^2$;
\item $\dim (U_n(x))\leqslant 2$.
\end{itemize}
\end{lema}

{\textbf{Proof}}. The first inequality is a consequence of 
Proposition~\ref{tec1}. 

Now, consider the nodal graph $\overline{\mathcal{N}}(f)$ in 
$\overline{\Sigma}$ of a function $f\in U_n(x)$ with highest $ord(f,x)$. 
Add new edges cutting non-simply-connected faces into discs (and hence 
making the graph connected). Contract some edges to merge all the vertices 
with $x$. We have obtained a new graph $\overline{\mathcal{N}}'$. If some 
of propositions of Lemma~\ref{cases} failed, namely, there were other 
vertices except $x$, or there were non-simply-connected faces 
of $\overline{\mathcal{N}}$, or $ord(f_x,x)>n$, or $\dim (U_n(x))>2$, 
then the degree of vertex $x$ in the new graph 
satisfies $\deg_{\overline{\mathcal{N}}'}(x)\geqslant 2n+2$ and 
hence $\overline{\mathcal{N}}'$ has at least $n+1$ edges. The number of 
faces of $\overline{\mathcal{N}}'$ is the same as for $\mathcal{N}$, 
no more than $n$. Now Euler characteristic of 
$\overline{\Sigma}\cong \mathbb{RP}^2$ can be estimated as 
$$1=\#\{\mbox{vertices of }{\overline{\mathcal{N}}'}\}-
\#\{\mbox{edges of }{\overline{\mathcal{N}}'}\}+
\#\{\mbox{faces of }{\overline{\mathcal{N}}'}\}\leqslant$$ $$\leqslant 1-
(n+1)+n=0,$$
but it is a contradiction.

The inequality $\dim (U(x))\leqslant 2$ follows now from the fact that, 
according to the proof of Proposition~\ref{tec1}, otherwise there would 
be a function $f_x\in U_{n+1}(x)$, such that $ord(f_x,x)>n$ which 
is already ruled out. 

Hence all the ways Lemma~\ref{cases} could fail, lead to the 
contradiction and Lemma~\ref{cases} is proven. $\square$

Lemma~\ref{cases} implies that if \mbox{$\dim (U)\geqslant 2n$} then for 
each point $x\in \mathbb{RP}^2$ there is either unique (up to 
$\mathbb{R}^*$) eigenfunction $f_x$ whose nodal graph $\mathcal{N}_{f_x}$ 
has the vertex $x$ of degree $2n$, or a 2-dimensional space $U(x)$ of 
such functions. We prove now that the case $\dim (U_n(x))=2$ described 
above is impossible in the case of odd $n$.

\begin{lema}
\label{1}
In the setting of Lemma~\ref{cases} if $n$ is an odd number then
we have got $\dim (U_n(x))=1$. 
\end{lema}

Before we get to the proof of Lemma~\ref{1}, let us quickly mention a 
technical property of nodal graphs in the surface $\overline{\Sigma}$ 
with the collapsed boundary.

\begin{prop}
\label{cont}
Suppose $f(p)$ is a continuous family of functions depending on a 
parameter $p$, each $f(p)$ satisfy $\Delta f(p)=\lambda f(p)$ for some 
real number $\lambda \in \mathbb{R}$ and the nodal 
graph $\overline{\mathcal{N}}(p)= \overline{\mathcal{N}}(f(p))$ has only
one vertex $x(p)$ in $\overline{\Sigma}$ 
and $\deg_{\overline{\mathcal{N}}(p)}(x(p))=2n$. Then the loops provided by 
edges of $\overline{\mathcal{N}}(p)$ do not change their homotopy class in 
the local system $\pi_1(\overline{\Sigma},x(p))$ while $p$ moves 
along $Int(\Sigma)$.
\end{prop}

The idea of the proof is that since $x(p_0)$ is the only vertex of the graph 
in $\overline{\Sigma}$, there is no more than two rays of the nodal graph 
in $\Sigma$ approaching each connected component of the boundary. Hence for 
all $p$ close to $p_0$ these rays can connect only with each other near 
the collapsed boundary component in $\overline{\Sigma}$. For details, 
see~\cite[Proposition~3]{berd}.

\textbf{Proof of Lemma~\ref{1}}.
Suppose that $\dim (U_n(x))=2$. Then, according to Bers theorem, 
for some polar coordinates $(r,\theta)$ near $x$ and polar 
coordinates $(R,\phi)$ in $U(x)$ we have the following approximation
$$
f_{(R,\phi)}(r,\theta)=Rr^{n}\cos((n)\theta+\phi)+O(r^{n}).
$$
Hence, taking the star of the function $\sigma(\bullet,x): PU\rightarrow F_x$ 
is a diffeomorphism of $PU\cong S^1\cong F_x$, and the generator loop of
$PU$ rotates each ray of $\sigma(f,x)$ to the subsequent one. 

Consider 
the universal cover 
$$\pi: \mathbb{R}\times [-1,1]\xrightarrow
{(\bullet /2\mathbb{Z})\times id} S^1\times[-1,1]\cong S^2\setminus 
\{\pm x\}\rightarrow$$
$$\xrightarrow{\bullet /\{\pm 1\}} \mathbb{R}P^2\setminus 
\{x\}\cong  S^1\ltimes [-1,1]$$ 
such that the fundamental group of $\mathbb{R}P^2\setminus \{x\}$ 
acts by integer shifts on $\mathbb{R}$-factor. Then the lift of the 
edge of $\overline{\mathcal{N}}(f)$ joins $(y,-1)$ and $(y+\frac{j}{n},e)$ 
for some $j\in\mathbb{Z}$ and $e\in\{-1,1\}$. Applying 
Proposition~\ref{cont} to an isotopy of a graph induced by a generator 
loop of $PU$, we conclude that all the points $(y+\frac{i}{n},-1)$ are 
joined with $(y+\frac{i+j}{n},e)$ by the lift of some edge 
of $\overline{\mathcal{N}}(f)$. Since this identification is a 
part of an involution prescribed by the lifts of the edges,  
$e$ has to be equal to 1. We conclude that $(y+\frac{i}{n},1)$ 
is joined with $(y+\frac{i-j}{n},-1)$. The lift of the antipodal 
map on $S^2\setminus \{\pm x\}$ replaces $j$ defined in this way 
by $-j$. But the graph, lifted along the quotient by $\{\pm 1\}$ 
is invariant under the antipodal map, therefore $j=0$. But this 
couldn't be either, as in this case for sufficiently small $\varepsilon$ we 
get 
$$
sgn(f(\pi(y+\frac{1}{2n},-(1-\varepsilon))))=
sgn(f(\pi(y+\frac{1}{2n},1-\varepsilon))=
$$
$$
=-sgn(f(\pi(y+\frac{n+1}{2n},1-\varepsilon)))=
-sgn(f(\pi(y+\frac{1}{2n},-(1-\varepsilon))))
$$
which is a contradiction. Here for the first equality we use $j=0$, for 
the second equality we use the fact that $n$ is odd and 
hence we have $sgn(\sin(\frac{\pi}{2}))=-sgn(\sin(n\pi+\sin(\frac{\pi}{2}))$ and 
for the third equality we use the invariance of $f$ under the action 
of $\pi_1(\mathbb{RP}^2 \setminus \{ x\} )$. $\square$

Finally, for an odd $n$ we are left with only possibility that for each 
point $x\in\mathbb{R}P^2$ there is unique (up to $\mathbb{R}^*$) 
eigenfunction $f_x$ such that $ord(f_x,x)=n$ and no functions 
with $ord(f_x,x)\geqslant n+1$. Note that it completes the proof of the 
first bound in the Theorem~\ref{B} required for the Theorem~\ref{m2}. 
Indeed, if in the conditions of the Theorem~\ref{B} we have 
$\dim (U)\geqslant 2n+1$, then by Proposition~\ref{tec1} we have
the inequality $\dim(U_{n}(x))\geqslant 1+(2n+1)-2n=2$, which 
contradicts Lemma~\ref{1} 
for odd $n$. 

Therefore Theorem~\ref{m2} is now proven. $\square$

\section{Second bound: closed surface case}

We switch to the case of $\Sigma=\mathbb{RP}^2$. Suppose that there 
is a counterexample to the Theorem~\ref{m1}, so that there is an 
eigenspace $U_{2l}$ with $\dim(U_{2l})\geqslant 4l+2$. Then due to 
Courant nodal domain theorem the assumption 
$\sup\limits_{f\in U}(\mu(\mathcal{N}(f)))\leqslant 2l+1=n$ of the 
Lemmas~\ref{cases} and~\ref{1} holds and for each point $x\in \mathbb{R}P^2$ 
we get a function $f_x$, defined up to a constant, such that $ord(f_x,x)=n$. 
The stars of the functions $f_x$ form then a smooth section $\sigma$ of
the star fibration $E_{\mathbb{RP}^2}(2n)$ according to the following 
technical proposition.

\begin{prop}
Let $U$ be a finite-dimensional eigenspace of the Laplacian. For every 
$x\in Int(\Sigma)$ consider the subspace $[f_x]\subset U$ of functions 
$f_x$ with $ord(f_x,x)\geqslant n$. Suppose that for any $x\in Int(\Sigma)$ 
the subspace $[f_x]$ is of dimension 1 (i. e. nonzero $f_x$ are defined 
up to $\mathbb{R}^*$), and that  $ord(f_x,x)=n$. Define a section 
$\sigma (x)=s(\mathcal{N}(f_x),x)$. Then this section 
$\sigma \in \Gamma (E_{Int(\Sigma)}(2n))$ is smooth.
\end{prop}

The proof follows directly from inverse function theorem, see 
the paper~\cite[Proposition 1]{berd}.


Now in order to achieve the final contradiction let us
pull $\sigma$ back to the star fibration $E_{\mathbb{S}^2}(2n)$ via the universal 
cover $\mathbb{S}^2\rightarrow \mathbb{RP}^2$ and observe that the Euler 
class $e(E_{\mathbb{S}^2}(2n))=2ne(T\mathbb{S}^2)=4n$ 
is non-zero and hence $E_{\mathbb{S}^2}(2n)$ can 
not have a continuous section. Hence, the initial assumption that the
multiplicity of $\lambda_{2l}$ is at least $4l+2$ is false and thus 
the multiplicity is 
at most $4k+1$.

Unfortunately, we did not find any evident way to make the analogous 
final consideration with the Euler class in the case 
$\partial \Sigma\neq \varnothing$, since the surface with boundary 
has no such class. However, it would be possible to use the relative 
Euler class if we could put some restrictions on the behavior of the 
star section $\sigma$ near the boundary. In the case of a surface $M$ of 
positive genus such restrictions come from the fact that some loops of 
nodal graph represent a (constant) non-zero class in $H_1(\tilde{M})$ 
for orienting cover $\tilde{M}$. This consideration, based on the method 
from the paper~\cite{closed}, is the crucial argument in
the paper~\cite{berd}. But in our case 
of $\mathbb{RP}^2$ the orienting cover $\mathbb{S}^2$ has 
zero homology in dimension 
1 and there is no apparent topological obstructions for the nodal graphs 
in the consideration to exist.

\end{document}